\documentclass{amsart}
\newtheorem{theorem}{Theorem}[section]
\newtheorem{corollary}[theorem]{Corollary}
\newtheorem{lemma}[theorem]{Lemma}
\newtheorem{prop}[theorem]{Proposition}
\newtheorem{conj}[theorem]{Conjecture}
\newcommand{\A}{\mathcal{A}}
\newcommand{\E}{\mathcal{E}}
\newcommand{\I}{\mathcal{I}}
\newcommand{\M}{\mathcal{M}}
\renewcommand{\O}{\mathcal{O}}
\renewcommand{\P}{\mathcal{P}}
\newcommand{\R}{\mathcal{R}}
\newcommand{\T}{\mathcal{T}}
\newcommand{\W}{\mathcal{W}}
\renewcommand{\b}{\underline{b}}
\renewcommand{\c}{\underline{c}}
\newcommand{\w}{\underline{w}}
\newcommand{\x}{\underline{x}}
\newcommand{\z}{\underline{z}}
\newcommand{\invlim}{\varprojlim}

\begin{document}
\title[Homeomorphisms of one-dimensional inverse limits]
{Homeomorphisms of One-dimensional Inverse
Limits with Applications to Substitution Tilings, Unstable Manifolds, and
Tent Maps}  

\author[M. Barge]{Marcy Barge}
\address{Department of Mathematical Sciences\\
Montana State University\\
Bozeman, Montana 59717}
\email{umsfbar@math.montana.edu}
\thanks{The first author was supported in part by NSF-DMS-9627026}

\author[J. Jacklitch]{James Jacklitch}
\email{jacklitc@math.montana.edu}
\thanks{The second author was supported in part by NSF-EPSCoR-93-50-546}

\author[G.M. Vago]{Gioia Vago}
\address{Laboratoire de Topologie\\ 
D\'epartement de Math\'ematiques\\
Universit\'e de Bourgogne\\
Dijon, France}
\email{vago@u-bourgogne.fr}
\thanks{The authors would like to thank Christian Bonatti for fruitful
conversations and the Laboratoire de Topologie, Universit\'e de Bourgogne, for
its warm hospitality.}


\begin{abstract}
Suppose that $f$ and $g$ are Markov surjections, each defined on a wedge
of circles, each fixing the branch point and having the branch
point as the only critical value. We show that if the points in the inverse 
limit spaces associated with $f$ and $g$ corresponding to the branch point are
distinguished then these inverse limit spaces are homeomorphic if and
only if the substitutions associated with $f$ and $g$ are weakly equivalent.
This, and related results, are applied to one-dimensional substitution tiling
spaces, one-dimensional unstable manifolds of hyperbolic sets, and 
inverse limits of tent maps with periodic critical points.
\end{abstract}

\maketitle
\setcounter{section}{-1}
\section{Introduction}
Many of the one-dimensional spaces that arise as invariant sets in dynamical
systems can be effectively modeled using inverse limits. This approach was
first employed by R.F. Williams in a description of one-dimensional
hyperbolic attractors \cite{W1}. In that characterization, Williams showed that 
some power of a diffeomorphism restricted to a one-dimensional hyperbolic 
attractor is
topologically conjugate to the natural (shift) homeomorphism on an inverse
limit whose bonding map is an immersion on a wedge of circles. Further
\cite{W2}, Williams
introduced the notion of shift equivalence of maps to classify such natural
homeomorphisms up to topological conjugacy. In this paper we address the
classification, up to homeomorphism, of inverse limit spaces whose bonding 
maps act on wedges of circles. A weaker notion than shift equivalence, 
namely weak equivalence, will arise in this context. As far as we know, the
relation of weak equivalence of matrices first arose in \cite{BD}. Consequences 
of weak equivalence of matrices and its relation to dimension groups 
were considered in
\cite{SV} and the notion has arisen in the classification of almost finite
dimensional algebras up to stable isomorphism (see \cite{BJKR}, where weak
equivalence of matrices, transposed, is called $C^*$-equivalence). Weak 
equivalence of matrices will play a part in the theorems we prove here but a
stronger equivalence relation, weak equivalence of substitutions, will play
the larger role. Our main result (Theorem \ref{T1.16}) asserts that certain inverse
limit spaces are homeomorphic if and only if the associated substitutions are
weakly equivalent.

Our motivation for considering inverse limits of maps on wedges of circles
comes from several sources. One is the classification of one-dimensional 
hyperbolic attractors, up to homeomorphism, as developed in \cite{J}. 
Others are: one-dimensional substitution tiling spaces;
one-dimensional unstable manifolds
of hyperbolic basic sets; and a conjecture of Ingram regarding inverse limits
of tent maps. These topics will be addressed in the final section of this
paper.
\section{Main Results}
For $n=1,2,\ldots,$ let $X_n=\bigvee_{i=1}^nS_i^1$ denote $n$ topological
circles, $S_i^1$, joined together at a common point, which we denote by $b$.
Assume that each $S_i^1$ is given an orientation. A map (by which we will
always mean a continuous function) $f:X\to Y$ is {\em monotone} provided
$f^{-1}(\{y\})$ is connected for each $y\in Y$ and a map
$h:X_n\to X_m$ is said to be 
{\em piecewise monotone} provided $X_n$ is a union of finitely many arcs
restricted to each of which $h$ is monotone.
A piecewise monotone map is {\em strictly piecewise monotone} if it is not
constant on any nontrivial arc. The point $c\in X_n$ is a {\em critical
point} of the map $h:X_n\to X_m$ provided $h$ is not one-to-one in
any neighborhood of $c$ and $v\in X_m$ is a {\em critical value} of $h$ if
$v=h(c)$ for some critical point $c\in X_n$.

For $n,m=1,2,\ldots,$ let $\M_{n,m}$ denote the collection of all piecewise
monotone maps $f:X_n\to X_m$ that satisfy: $f(b)=b$; $b$ is the
only critical value (if any) of $f$; $f$ is surjective; and $f$ is not
constant on any $S_i^1$, $i=1,\ldots,n$. Let
$\M=\bigcup_{n=1}^{\infty}\M_{n,n}$. We are interested in the 
topology of the {\em inverse limit space}, $\invlim f$, with {\em bonding maps}
$f\in\M$. More 
generally, if $f_i:X_n\to X_n$, $i=1,2,\ldots,$ are maps in $\M_{n,n}$, we will
denote by $\invlim f_i$ the space $\{(x_0,x_1,\ldots)\in\prod_{k=0}^{\infty}X_n:
f_i(x_i)=x_{i-1}\;{\rm for}\; i=1,2,\ldots\}$. The
topology on $\invlim f_i$ is inherited from the product space
$\prod_{k=0}^{\infty}X_n$; with this topology $\invlim f_i$ is a {\em continuum}
(compact, connected, metrizable). We will denote by $\pi_k$ the projection
$\pi_k:\invlim f_i\to X_n$, $\pi_k(x_0,x_1,\ldots,x_k,\ldots)=x_k$
and, when all bonding maps are the same, $f_i=f$ for each $i$,
$\hat{f}:\invlim f\to \invlim f$ will denote the natural (shift)
homeomorphism $\hat{f}(x_0,x_1,\ldots)=(f(x_0),x_0,x_1,\ldots)$.

For each $n=1,2,\ldots,$ let $\A_n$ denote the alphabet consisting of $2n$
letters $\A_n=\{a_1,\ldots,a_n,\bar{a}_1,\ldots,\bar{a}_n\}$ and let
$\W(\A_n)$ be the collection of all finite nonempty words with letters in
$\A_n$. Given the $k$-letter word $b_1\ldots b_k\in\W(\A_m)$,
let $\overline{b_1\ldots b_k}=\bar{b}_k\ldots \bar{b}_1\in\W(\A_m)$
with the convention that ${\bar{\bar{b}}}_i=b_i$. 
Starting from $f\in\M_{n,m}$, we want to define an associated substitution
$\chi_f:\A_n\to\W(\A_m)$ in such a way that $\chi_f(a_i)$ (resp.
$\chi_f(\bar{a}_i)$) describes the way in which $f$ maps the circle $S_i^1$
(resp. the circle $S_i^1$ oriented in the opposite way) around each $S_j^1$,
keeping track of order and orientation. For this reason, denote by 
$\dot{J}_1<\dot{J}_2<\ldots<\dot{J}_k$ the components of $S_i^1\backslash
f^{-1}(\{b\})$, the ordering of such arcs determined by the positive
orientation on $S_i^1$. Then define  $\chi_f(a_i)=b_1\ldots b_k$ provided $f$
maps $\dot{J}_j$ homeomorphically onto $\dot{S}_{\ell}^1\equiv
S_{\ell}^1\backslash\{b\}$ in an orientation preserving (resp. reversing) way
if and only if $b_j=a_{\ell}$ (resp. $\bar{a}_{\ell}$). Recall the meaning of
a substitution. Now, when $S_i^1$ is oriented according to the opposite
orientation, $S_i^1\backslash f^{-1}(\{b\})$ has components $\dot{J}_k
<\dot{J}_{k-1}<\ldots<\dot{J}_1$. Then it is natural to define 
$\chi_f(\bar{a}_i)\equiv \bar{b}_k\ldots\bar{b}_1=\overline{\chi_f(a_i)}$.
We then extend $\chi_f$ to $\chi_f:\W(\A_n)\to\W(\A_m)$ by
concatenation: $\chi_f(c_1\ldots c_{\ell})=\chi_f(c_1)
\ldots\chi_f(c_{\ell})$. Note that $\chi_f(\overline{c_1\ldots c_{\ell}})=
\overline{\chi_f(c_1\ldots c_{\ell})}$. For any substitution $\tau:\A_n
\to\W(\A_m)$ satisfying $\tau(\bar{b})=
\overline{\tau(b)}$, let $A_{\tau}$ be the $m\times n$ matrix whose $ji$-th
entry is $(A_{\tau})_{ji}=\#\{\ell:b_{\ell}\in \{a_j,\bar{a}_j\}\;{\rm where}\;
\tau(a_i)=b_1\ldots b_k\}$. That is, $(A_{\tau})_{ji}$ is the number of
occurrences of the letters $a_j$ and $\bar{a}_j$ in the word $\tau(a_i)$. In
case $\tau=\chi_f$, we will say that $f$ {\em follows the pattern} $\tau$ and we
will write $A_f$ for $A_{\tau}$.

Suppose that $\chi:\A_n\to\W(\A_n)$ and $\psi:\A_m\to\W(\A_m)$ 
are substitutions that
satisfy $\chi(\bar{b})=\overline{\chi(b)}$ and
$\psi(\bar{c})=\overline{\psi(c)}$ for all $b\in\A_n$, $c\in\A_m$.
We will say that $\chi$ and $\psi$ are {\em weakly equivalent}, $\chi\sim_w
\psi$, if there are sequences of positive integers $\{n_i\}$, $\{m_i\}$, and
substitutions $\sigma_i:\A_m\to\W(\A_n)$,
$\tau_i:\A_n\to\W(\A_m)$, $i=1,2,\ldots,$ such
that $\sigma_i(\bar{c})=\overline{\sigma_i(c)}$, $\tau_i(\bar{b})=
\overline{\tau_i(b)}$ for all $c\in\A_m$, $b\in\A_n$ and, after extending
the substitutions by concatenation,
$\sigma_i \circ \tau_i=\chi^{n_i}$, $\tau_i \circ \sigma_{i+1}=\psi^{m_i}$.
That is, the following infinite diagram commutes.
\begin{center}
\begin{picture}(4.8,1.8)
\put(1.35,.4){\makebox(0,0){$\W(\A_m)$}}
\put(.4,1.4){\makebox(0,0){$\W(\A_n)$}}
\put(3.25,.4){\makebox(0,0){$\W(\A_m)$}}
\put(2.3,1.4){\makebox(0,0){$\W(\A_n)$}}
\put(.7,.8){\makebox(0,0){$\sigma_1$}}
\put(1.6,1){\makebox(0,0){$\tau_1$}}
\put(2.6,.8){\makebox(0,0){$\sigma_2$}}
\put(3.5,1){\makebox(0,0){$\tau_2$}}
\put(2.4,.2){\makebox(0,0){$\psi^{m_1}$}}
\put(1.45,1.6){\makebox(0,0){$\chi^{n_1}$}}
\put(4.3,.2){\makebox(0,0){$\psi^{m_2}$}}
\put(3.35,1.6){\makebox(0,0){$\chi^{n_2}$}}
\put(4.7,.9){\makebox(0,0){$\ldots$}}
\put(1.05,.6){\vector(-1,2){.3}}
\put(2,1.2){\vector(-1,-2){.3}}
\put(2.95,.6){\vector(-1,2){.3}}
\put(3.9,1.2){\vector(-1,-2){.3}}
\put(2.775,.4){\vector(-1,0){.95}}
\put(1.85,1.4){\vector(-1,0){1}}
\put(4.675,.4){\vector(-1,0){.95}}
\put(3.75,1.4){\vector(-1,0){1}}
\end{picture}
\end{center}

If $A_{n\times n}$ and $B_{m\times m}$ are square matrices with nonnegative
integer entries, we will say that $A$ and $B$ are {\em weakly equivalent}, $A\sim_w
B$, provided there are sequences of positive integers $\{n_i\}$, $\{m_i\}$,
and nonnegative integer matrices $S_i$, $T_i$, $i=1,2,\ldots$, such that
$S_iT_i=A^{n_i}$ and $T_iS_{i+1}=B^{m_i}$. (Weak equivalence of $A$ and $B$
is called $C^*$-equivalence of $A$ transpose and $B$ transpose in \cite{BJKR}.)

\begin{lemma}\label{L1.1}
If $f,g\in\M$ and $\chi_f\sim_w \chi_g$ then $A_f\sim_w A_g$.
\end{lemma}

\begin{proof} 
If ${n_i}$, ${m_i}$, $\sigma_i$, $\tau_i$ realize the weak
equivalence of $\chi_f$ and $\chi_g$, let $S_i$ and $T_i$ be given by
$S_i=A_{\sigma_i}$ and $T_i=A_{\tau_i}$. It is immediate that $A_f^{n_i}=S_iT_i$
and $A_g^{m_i}=T_iS_{i+1}$.
\end{proof}

\begin{lemma}\label{L1.2}
Given $f\in\M_{n,n}$ there are strictly piecewise monotone
$f_i\in\M_{n,n}$, $i=1,2,\ldots,$ each of which follows the pattern of
$\chi_f$, such that $\invlim f_i\simeq\invlim f$.
\end{lemma}

\begin{proof}
The approximation theorem of Mort Brown (\cite{Bro}) asserts that there
are $\epsilon_i>0$ so that if $f_i:X_n\to X_n$ is uniformly within
$\epsilon_i$ of $f$ (fix some metric on $X_n$) for each $i$ then
$\invlim f_i \simeq \invlim f$. Clearly, such $f_i$ can be
chosen to be strictly piecewise monotone, in $\M_{n,n}$, and follow the
pattern of $\chi_f$.
\end{proof} 

\begin{lemma}\label{L1.3}
If $f,g\in\M_{n,n}$ and $\chi_f=\chi_g$ then $\invlim f\simeq\invlim g$. 
\end{lemma}

\begin{proof}
Let $f_i,g_i\in\M_{n,n}$ be as in Lemma \ref{L1.2} with
$\invlim f_i \simeq \invlim f$ and $\invlim g_i \simeq
\invlim g$, all following the pattern of $\chi_f=\chi_g$. Let
$P_0:X_n\to X_n$ be the identity. Suppose that, for $j=0,\ldots,k$,
homeomorphisms $P_j:X_n\to X_n$ have been defined in such a way that 
$P_j(S_i^1)=S_i^1$ for each $i\in\{1,\ldots,n\}$ and the diagram
\begin{center}
\begin{picture}(5.6,2.2)
\put(.4,.4){\makebox(0,0){$X_n$}}
\put(.4,1.8){\makebox(0,0){$X_n$}}
\put(1.8,.4){\makebox(0,0){$X_n$}}
\put(1.8,1.8){\makebox(0,0){$X_n$}}
\put(5.2,.4){\makebox(0,0){$X_n$}}
\put(5.2,1.8){\makebox(0,0){$X_n$}}
\put(.2,1.1){\makebox(0,0){$P_0$}}
\put(1.6,1.1){\makebox(0,0){$P_1$}}
\put(5,1.1){\makebox(0,0){$P_k$}}
\put(1.1,.2){\makebox(0,0){$g_1$}}
\put(1.1,2){\makebox(0,0){$f_1$}}
\put(2.5,.2){\makebox(0,0){$g_2$}} 
\put(2.5,2){\makebox(0,0){$f_2$}}
\put(4.5,.2){\makebox(0,0){$g_k$}}
\put(4.5,2){\makebox(0,0){$f_k$}}
\put(3.5,1.1){\makebox(0,0){$\ldots$}}
\put(.4,1.6){\vector(0,-1){1}}
\put(1.8,1.6){\vector(0,-1){1}}
\put(5.2,1.6){\vector(0,-1){1}}
\put(1.6,.4){\vector(-1,0){1}}
\put(1.6,1.8){\vector(-1,0){1}}
\put(3,.4){\vector(-1,0){1}}
\put(3,1.8){\vector(-1,0){1}}
\put(5,.4){\vector(-1,0){1}}
\put(5,1.8){\vector(-1,0){1}}
\end{picture}
\end{center}
commutes. Since $f_{k+1},g_{k+1}\in\M_{n,n}$ are strictly piecewise
monotone, and they follow the same pattern, for each $i$, $S_i^1\backslash
f^{-1}_{k+1}(\{b\})=\bigcup_{\ell=1}^{r_i}\dot{J}_{i,\ell}$ and $S_i^1\backslash
g^{-1}_{k+1}(\{b\})=\bigcup_{\ell=1}^{r_i}\dot{J}_{i,\ell}^{\prime}$ where the 
$\dot{J}_{i,\ell}$ and $\dot{J}_{i,\ell}^{\prime}$ are pairwise disjoint 
open arcs, $\bigcup_{\ell=1}^{r_i}\dot{J}_{i,\ell}$ and 
$\bigcup_{\ell=1}^{r_i}\dot{J}_{i,\ell}^{\prime}$
are dense in $S_i^1$, $\dot{J}_{i,1}<\dot{J}_{i,2}<\ldots<\dot{J}_{i,r_i}$
and $\dot{J}_{i,1}^{\prime}<\dot{J}_{i,2}^{\prime}<\ldots<
\dot{J}_{i,r_i}^{\prime}$ (the order given by the positive orientation on
$S_i^1$), $f_{k+1}|_{\dot{J}_{i,\ell}}$ and
$g_{k+1}|_{\dot{J}_{i,\ell}^{\prime}}$ are both orientation preserving 
or both orientation reversing homeomorphisms, and
$f_{k+1}(\dot{J}_{i,\ell})=\dot{S}_r^1=g_{k+1}(\dot{J}_{i,\ell}^{\prime})$ 
for some $r=r(i,\ell)$. Let $P_{k+1}|_{\dot{J}_{i,\ell}}:\dot{J}_{i,\ell}\to
\dot{J}_{i,\ell}^{\prime}$ be defined by
$P_{k+1}|_{\dot{J}_{i,\ell}}=(g_{k+1}|_{\dot{J}_{i,\ell}^{\prime}})^{-1} \circ P_k 
\circ f_{k+1}|_{\dot{J}_{i,\ell}}$. $P_{k+1}$ then extends from
$\bigcup_{i=1}^n\bigcup_{j=1}^{r_i}\dot{J}_{i,j}$ to a homeomorphism of $X_n$
onto itself so that $P_k\circ f_{k+1}=g_{k+1}\circ P_{k+1}$. 
\begin{center}
\begin{picture}(7.6,2.2)
\put(.4,.4){\makebox(0,0){$X_n$}}
\put(.4,1.8){\makebox(0,0){$X_n$}}
\put(1.8,.4){\makebox(0,0){$X_n$}}
\put(1.8,1.8){\makebox(0,0){$X_n$}}
\put(5.2,.4){\makebox(0,0){$X_n$}}
\put(5.2,1.8){\makebox(0,0){$X_n$}}
\put(.2,1.1){\makebox(0,0){$P_0$}}
\put(1.6,1.1){\makebox(0,0){$P_1$}}
\put(5,1.1){\makebox(0,0){$P_k$}}
\put(1.1,.2){\makebox(0,0){$g_1$}}
\put(1.1,2){\makebox(0,0){$f_1$}}
\put(2.5,.2){\makebox(0,0){$g_2$}}
\put(2.5,2){\makebox(0,0){$f_2$}}
\put(4.5,.2){\makebox(0,0){$g_k$}}
\put(4.5,2){\makebox(0,0){$f_k$}}
\put(5.9,.2){\makebox(0,0){$g_{k+1}$}}
\put(5.9,2){\makebox(0,0){$f_{k+1}$}}
\put(3.5,1.1){\makebox(0,0){$\ldots$}}
\put(6.9,1.1){\makebox(0,0){$\ldots$}}
\put(.4,1.6){\vector(0,-1){1}}
\put(1.8,1.6){\vector(0,-1){1}}
\put(5.2,1.6){\vector(0,-1){1}}
\put(1.6,.4){\vector(-1,0){1}}
\put(1.6,1.8){\vector(-1,0){1}}
\put(3,.4){\vector(-1,0){1}}
\put(3,1.8){\vector(-1,0){1}}
\put(5,.4){\vector(-1,0){1}}
\put(5,1.8){\vector(-1,0){1}}
\put(6.4,.4){\vector(-1,0){1}}
\put(6.4,1.8){\vector(-1,0){1}}
\end{picture}
\end{center}
The preceding commuting diagram
then induces a homeomorphism $\hat{P}:\invlim f_i\to\invlim g_i$ by
$\hat{P}((x_0,x_1,\ldots))=(P_0(x_0),P_1(x_1),\ldots)$. 
\end{proof}

\begin{lemma}\label{L1.4}
Given $f\in\M_{n,n}$ there is a strictly piecewise monotone
$f^{\prime}\in\M_{n,n}$ that follows the same pattern as $f$ with $\invlim 
f^{\prime}\simeq\invlim f$.
\end{lemma}

\begin{proof}
This follows immediately from Lemma \ref{L1.3} by taking, for instance,
$f^{\prime}$ to be any of the $f_i$ of \ref{L1.2}.
\end{proof}

\begin{theorem}\label{T1.5}
If $f,g\in\M$ with $\chi_f\sim_w\chi_g$ then $\invlim f\simeq\invlim g$.
\end{theorem}

\begin{proof}
By Lemma \ref{L1.4}, we may assume that $f$ and $g$ are strictly 
piecewise monotone. Suppose that $f\in\M_{n,n}$, $g\in\M_{m,m}$ and let $\{n_i\}$,
$\{m_i\}$, $\sigma_i:\A_m\to\W(\A_n)$, $\tau_i:\A_n\to\W(\A_m)$ be such that,
after extending the substitutions by concatenation,
$\sigma_i\tau_i=\chi_f^{n_i}$, $\tau_i\sigma_{i+1}=\chi_g^{m_i}$. We will
construct continuous surjections $s_i:X_m\to X_n$ and $t_i:X_n\to X_m$ such
that $s_i\circ t_i=f^{n_i}$ and $t_i\circ s_{i+1}=g^{m_i}$.

Let $s_1:X_m\to X_n$ be any continuous strictly piecewise monotone surjection
that takes $b$ to $b$ and follows the pattern $\sigma_1$. Suppose that
$s_i:X_m\to X_n$, $i=1,\ldots,N$, and $t_i:X_n\to X_m$, $i=1,\ldots,N-1$, are
all continuous strictly piecewise monotone surjections taking $b$ to $b$ following
the patterns $\sigma_i$ and $\tau_i$ that have been defined so that the 
following diagram commutes.
\begin{center}
\begin{picture}(6.3,1.8)
\put(1.1,.4){\makebox(0,0){$X_m$}}
\put(.4,1.4){\makebox(0,0){$X_n$}}
\put(2.5,.4){\makebox(0,0){$X_m$}}
\put(1.8,1.4){\makebox(0,0){$X_n$}}
\put(5.9,.4){\makebox(0,0){$X_m$}}
\put(5.2,1.4){\makebox(0,0){$X_n$}}
\put(.6,.8){\makebox(0,0){$s_1$}}
\put(1.25,1){\makebox(0,0){$t_1$}}
\put(2,.8){\makebox(0,0){$s_2$}}
\put(2.65,1){\makebox(0,0){$t_2$}}
\put(4.45,1){\makebox(0,0){$t_{N-1}$}}
\put(5.3,.8){\makebox(0,0){$s_N$}}
\put(1.9,.2){\makebox(0,0){$g^{m_1}$}}
\put(1.2,1.6){\makebox(0,0){$f^{n_1}$}}
\put(3.3,.2){\makebox(0,0){$g^{m_2}$}}
\put(2.6,1.6){\makebox(0,0){$f^{n_2}$}}
\put(5.3,.2){\makebox(0,0){$g^{m_{N-1}}$}}
\put(4.6,1.6){\makebox(0,0){$f^{n_{N-1}}$}}
\put(3.7,.9){\makebox(0,0){$\ldots$}}
\put(.9,.6){\vector(-1,2){.3}}
\put(1.6,1.2){\vector(-1,-2){.3}}
\put(2.3,.6){\vector(-1,2){.3}}
\put(3,1.2){\vector(-1,-2){.3}}
\put(5,1.2){\vector(-1,-2){.3}}
\put(5.7,.6){\vector(-1,2){.3}}
\put(2.3,.4){\vector(-1,0){1}}
\put(1.6,1.4){\vector(-1,0){1}}
\put(3.7,.4){\vector(-1,0){1}}
\put(3,1.4){\vector(-1,0){1}}
\put(5.7,.4){\vector(-1,0){1}}
\put(5,1.4){\vector(-1,0){1}}
\end{picture}
\end{center}

Now $(f^{n_N})^{-1}(\{b\})$ breaks each $S_i^1\subset X_n$ into a certain
number of intervals; we define the same number of intervals in $X_m$ in two
steps, with the aid of the substitutions $\sigma_N$ and $\tau_N$.
Fix $i\in\{1,\ldots,n\}$ and suppose that $\chi_f^{n_N}(a_i)=b_1\ldots b_s$.
Then $(f^{n_N})^{-1}(\{b\})\cap S_i^1=\{b=x_0<x_1<\ldots<x_s=b\}$ with
$f^{n_N}$ mapping $(x_{j-1},x_j)$ homeomorphically onto $\dot{S}_{\ell}^1$ in an
orientation preserving (resp. reversing) way iff 
$b_j=a_{\ell}$ (resp. $\bar{a}_{\ell}$). Let
$\tau_N(a_i)=c_1\ldots c_k$. Since
$\sigma_N\tau_N(a_i)=\chi_f^{n_N}(a_i)=b_1\ldots b_s$, there are integers
$0=p_0<p_1<\ldots<p_k=s$ such that $\sigma_N(c_r)=b_{p_{r-1}+1}\ldots b_{p_r}$
for $r=1,\ldots,k$. Now fix $j\in\{1,\ldots,s\}$ and let $r$ be such that
$p_{r-1}<j\le p_r$. Suppose that $c_r\in\{a_q,\bar{a}_q\}$. Since $s_N$
follows the pattern $\sigma_N$, $s_N^{-1}(\{b\})\cap
S_q^1=\{b=y_0<y_1<\ldots<y_{p_r-p_{r-1}}=b\}$ and $s_N$ maps
$(y_{j-p_{r-1}-1},y_{j-p_{r-1}})$ homeomorphically onto $\dot{S}_{\ell}^1$ 
where $b_j\in\{a_{\ell},\bar{a}_{\ell}\}$. Define $t_N$ on $(x_{j-1},x_j)$ 
by $t_N=(s_N|_{(y_{j-p_{r-1}-1},y_{j-p_{r-1}})})^{-1}\circ
f^{n_N}|_{(x_{j-1},x_j)}$. Extend $t_N$ to a continuous surjection from
$[x_{j-1},x_j]$ onto $[y_{j-p_{r-1}-1},y_{j-p_{r-1}}]$. Carry out the above
for each $i\in\{1,\ldots,n\}$ and $j\in\{1,\ldots,s=s_i\}$. That there is no
ambiguity in the definition of $t_N$ is guaranteed by the construction and
the fact that
$\sigma_N\tau_N=\chi_f^{n_N}$. It is clear that $t_N$ is a
continuous strictly piecewise monotone surjection that follows the pattern
$\tau_N$ and that $s_Nt_N=f^{n_N}$. Now define $s_{N+1}:X_m\to X_n$ using
$t_N$ and $g^{m_N}$, etc. The commuting diagram
\begin{center}
\begin{picture}(4,1.8)
\put(1.1,.4){\makebox(0,0){$X_m$}}
\put(.4,1.4){\makebox(0,0){$X_n$}}
\put(2.5,.4){\makebox(0,0){$X_m$}}
\put(1.8,1.4){\makebox(0,0){$X_n$}}
\put(.6,.8){\makebox(0,0){$s_1$}}
\put(1.25,1){\makebox(0,0){$t_1$}}
\put(2,.8){\makebox(0,0){$s_2$}}
\put(2.65,1){\makebox(0,0){$t_2$}}
\put(1.9,.2){\makebox(0,0){$g^{m_1}$}}
\put(1.2,1.6){\makebox(0,0){$f^{n_1}$}}
\put(3.3,.2){\makebox(0,0){$g^{m_2}$}}
\put(2.6,1.6){\makebox(0,0){$f^{n_2}$}}
\put(3.7,.9){\makebox(0,0){$\ldots$}}
\put(.9,.6){\vector(-1,2){.3}}
\put(1.6,1.2){\vector(-1,-2){.3}}
\put(2.3,.6){\vector(-1,2){.3}}
\put(3,1.2){\vector(-1,-2){.3}}
\put(2.3,.4){\vector(-1,0){1}}
\put(1.6,1.4){\vector(-1,0){1}}
\put(3.7,.4){\vector(-1,0){1}}
\put(3,1.4){\vector(-1,0){1}}
\end{picture}
\end{center}
induces a homeomorphism between $\invlim f$ and $\invlim g$.
\end{proof}

In order to obtain a converse to Theorem \ref{T1.5} we will need to assume 
that the point $\b=(b,b,\ldots)$ in $\invlim f$ and $\invlim g$ is somehow
distinguished so that any homeomorphism would
necessarily take $\b$ to $\b$. We proceed to isolate three ways this can happen
for $f,g\in\M$. To this end, consider the action of $f\in\M_{n,n}$ on the ``edge
germs" at $b$. There are $2n$ germs at $b$, 
$\E=\{e_1,\ldots,e_n,\bar{e}_1,\ldots,\bar{e}_n\}$, where we denote the edge
germs of $S_i^1$ at $b$ by $e_i$ and $\bar{e}_i$: $e_i$ (resp. $\bar{e}_i$)
emanates from $b$ in the direction of the positive (resp. negative)
orientation of $S_i^1$. Define $f_*:\E\to\E$ by $f_*(e_i)=e_j$ (resp.
$\bar{e}_j$) provided $\chi_f(a_i)=a_j\ldots$ (resp. $\bar{a}_j\ldots$) and
$f_*(\bar{e}_i)=e_j$ (resp. $\bar{e}_j$) provided
$\chi_f(\bar{a}_i)=a_j\ldots$ (resp. $\bar{a}_j\ldots$).
Let $\R(f)=\bigcap_{k\ge 0}f_*^k(\E)$ be the eventual
range of $f_*$. Note that if $f\in\M_{n,n}$ then $\R(f)=f_*^{2n}(\E)$.

By an {\em open arc} we will mean a homeomorphic copy of the interval
$(-1,1)$.

\begin{lemma}\label{L1.6}
If $f\in\M$ and $\#\R(f)=1$ then there does not exist an open
arc in $\invlim f$ containing $\b$.
\end{lemma}

\begin{proof}
By Lemma \ref{L1.4}, we may assume that $f\in\M_{n,n}$ is strictly
piecewise monotone (the homeomorphisms of Lemmas \ref{L1.2} and \ref{L1.3}
can be chosen to take $\b$ to $\b$). Moreover, we may assume that $\#f_*(\E)=1$, 
since $\invlim f^{2n}\simeq\invlim f$, and that $f|_{[b,\epsilon]}$ is the
identity for some $\epsilon\ne b$ where $[b,\epsilon]$ is a closed arc in
$S_j^1$ and $f_*(\E)=\{e_j\}$ or $\{\bar{e}_j\}$. Then, if $\alpha$ is an open
arc in $\invlim f$ containing $\b$, we may assume that
$\alpha_k=\pi_k(\alpha)\subset[b,\epsilon]$ for all $k$. But then
$\alpha\subset\invlim f|_{[b,\epsilon]}$, and this is clearly not possible
since $\b$ is an endpoint of the closed arc $\invlim f|_{[b,\epsilon]}$.
\end{proof}

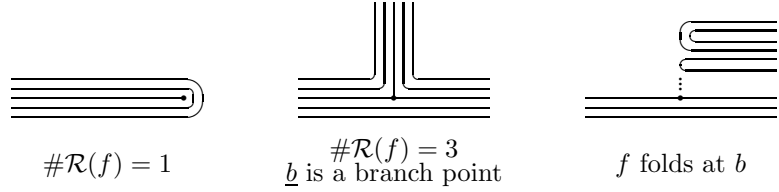
\begin{center}
\begin{picture}(8,2)
\put(1.8,1){\circle*{.06}}
\put(0,1){\line(1,0){1.8}}
\put(0,1){\oval(3.8,.2)[r]}
\put(0,1){\oval(4,.4)[r]}
\put(3,1){\line(1,0){2}}
\put(3,.9){\line(1,0){2}}
\put(3,.8){\line(1,0){2}}
\put(1,.3){\makebox(0,0){$\#\R(f)=1$}}
\put(4,.45){\makebox(0,0){$\#\R(f)=3$}}
\put(4,.2){\makebox(0,0){$\b$ is a branch point}}
\put(7,.3){\makebox(0,0){$f$ folds at $b$}}
\put(4,1){\circle*{.06}}
\put(4,1){\line(0,1){1}}
\put(3.8,1.1){\line(-1,0){.8}}
\put(3.7,1.2){\line(-1,0){.7}}
\put(3.9,1.2){\line(0,1){.8}}
\put(3.8,1.3){\line(0,1){.7}}
\put(4.2,1.1){\line(1,0){.8}}
\put(4.3,1.2){\line(1,0){.7}}
\put(4.1,1.2){\line(0,1){.8}}
\put(4.2,1.3){\line(0,1){.7}}
\put(3.8,1.2){\oval(.2,.2)[br]}
\put(3.7,1.3){\oval(.2,.2)[br]}
\put(4.2,1.2){\oval(.2,.2)[bl]}
\put(4.3,1.3){\oval(.2,.2)[bl]}
\put(7,1){\circle*{.06}}
\put(7,1.2){\circle*{.02}}
\put(7,1.15){\circle*{.02}}
\put(7,1.1){\circle*{.02}}
\put(6,1){\line(1,0){2}}
\put(6,.9){\line(1,0){2}}
\put(6,.8){\line(1,0){2}}
\put(8,1.35){\oval(2,.12)[l]}
\put(8,1.65){\oval(2,.3)[l]}
\put(8,1.65){\oval(1.8,.12)[l]}
\end{picture}
\end{center}

We will call a point $x$ in a continuum $X$ a {\em branch point} if there is
a continuum $Y\subset X$ such that $Y$ is homeomorphic with the letter Y
(three closed arcs joined at an endpoint) and $Y\backslash\{x\}$ has three
components.

\begin{lemma}\label{L1.7}
If $f\in\M$ and $\#\R(f)\ge 3$ then $\b$ is a branch point of $\invlim f$.
\end{lemma}

\begin{proof}
Replacing $f$ by $f^{2n!}$, if necessary, we may suppose that
there are distinct edge germs $b_i\in\E$ such that $f_*(b_i)=b_i$, for
$i=1,2,3$. Say, without loss of generality, that $b_i=e_i$, $i=1,2,3$. We may
also assume that $\epsilon_i\in\dot{S}_i^1$ are such that
$f|_{[b,\epsilon_i]}$ is the identity for $i=1,2,3$. Let $Y=\invlim
f|_{\cup_{i=1}^3[b,\epsilon_i]}$.
\end{proof}

We will say $f\in\M_{n,n}$ {\em folds at b} if there are
$i,j\in\{1,\ldots,n\}$ and $k\ge 1$ such that the letters $a_j$ and
$\bar{a}_j$ occur consecutively in $\chi_f^k(a_i)$.

\begin{lemma}\label{L1.8}
If $f\in\M_{n,n}$ folds at $b$ then no neighborhood of $\b\in\invlim f$ is
homeomorphic with the product of a zero-dimensional space and an open arc.
\end{lemma}

\begin{proof}
We may assume that $f_*$ fixes each edge germ in $\R(f)$ (by
replacing $f$ by $f^{2n!}$). For the sake of argument, assume that
$f_*(e_1)=e_1$ and that $\bar{a}_1$, $a_1$ occur consecutively in
$\chi_f(a_j)$. Further, assume that $f$ is strictly piecewise monotone and
that $f|_{[b,\epsilon]}$ is the identity for some $\epsilon\in
\dot{S}_1^1$. Let $a^{\prime}<c<b^{\prime}\in\dot{S}_j^1$ be such
that $f(c)=b$, $f(a^{\prime})=f(b^{\prime})=\epsilon$,
$f([a^{\prime},b^{\prime}])=[b,\epsilon]$. Let $J_0=[b,\epsilon]$,
$J_1=[a^{\prime},b^{\prime}]$, and choose intervals $J_2,J_3,\ldots$ such
that $f$ takes $J_{i+1}$ homeomorphically onto $J_i$ for $i\ge 1$. Let
$X\equiv\{(x_0,x_1,\ldots)\in\invlim f : {\rm for\; some}\;n, 
x_0=x_1=\ldots=x_n\;{\rm and}\;x_{n+k}\in J_k\;{\rm for}\;k\ge 1\}$.
Then every neighborhood of $\b$ in $\invlim f$ contains a
neighborhood of $\b$ in $X$, which in turn
contains a homeomorphic copy of $X$. Note that $X$ is homeomorphic with a
neighborhood of the point $(0,1)$ in the ``topologist's sine curve"
$\{0\}\times[-1,1]\cup\{(t,sin 1/t):0<t\le 1\}$. No such $X$ is contained
in the product of a zero-dimensional set with an arc.
\end{proof} 

Following \cite{AM} we will call a closed subset $M$ of a continuum $X$ a 
{\em matchbox in} $X$ if
there exist a zero-dimensional space $S$ and a topological embedding
$h:S\times[-1,1]\to X$ such that $h(S\times[-1,1])=M$ and
$h(S\times(-1,1))=\dot{M}$, the interior of $M$. 
The map $h:S\times[-1,1]\to M$ is called a {\em
parametrization of} $M$. An {\em oriented matchbox} is a matchbox provided
with a parametrization, whose matches (arc components) are oriented by the
canonical orientation of $[-1,1]$ via the parametrization. Let
$pr_1:S\times[-1,1]\to S$ be defined by $pr_1(s,t)=s$ and $pr_2:S\times[-1,1]\to
[-1,1]$ be defined by $pr_2(s,t)=t$. If $M_1$ and $M_2$ are oriented matchboxes 
in $X$ with parametrizations $h_1:S_1\times[-1,1]\to M_1$ and
$h_2:S_2\times[-1,1]\to M_2$ then we say that $M_1$ and $M_2$
are {\em oriented coherently} (resp. {\em anticoherently}) if for
each $x\in int(M_1\cap M_2)$ the map $pr_2|_{h_2^{-1}(A_x)}\circ
h_2^{-1}\circ h_1\circ (pr_2|_{h_1^{-1}(A_x)})^{-1}:pr_2\circ
h_1^{-1}(A_x)\to pr_2\circ h_2^{-1}(A_x)$ is increasing (resp. decreasing),
where $A_x$ is the component of $M_1\cap M_2$ containing $x$.

\begin{lemma}\label{L1.9}
If $f\in\M_{n,n}$ and $J$ is a closed arc of $X_n\backslash\{b\}$
then for $k=0,1,\ldots,$ $\pi_k^{-1}(J)$ is an oriented matchbox in $\invlim
f\backslash\{\b\}$.
\end{lemma}

\begin{proof}
Let $h:[-1,1]\to J$ be an orientation preserving homeomorphism onto $J$ (the
orientation of $J$ agreeing with the orientation of the $S_i^1$ containing it).  
Let $\dot{J}_1,\ldots,\dot{J}_{\ell}$ be the arc components of $X_n\backslash 
f^{-1}(\{b\})$ and $\Sigma_{\ell}$ be the product space
$\prod_{j=1}^{\infty}\{1,\ldots,\ell\}$ where $\{1,\ldots,\ell\}$ has the
discrete topology. Then $S\equiv\{(b_1,b_2,\ldots)\in\Sigma_{\ell}:J\subset 
f(\dot{J}_{b_1})\;{\rm and}\; \dot{J}_{b_i}\subset f(\dot{J}_{b_{i+1}})\;{\rm for}\;
i=1,2,\ldots\}$ is a closed subset of the Cantor set $\Sigma_{\ell}$. Let
$h_k:S\times[-1,1]\to\pi_k^{-1}(J)$ be the homeomorphism defined by 
 $h_k((b_1,b_2,\ldots),t)=(f^k\circ
h(t),f^{k-1}\circ h(t),\ldots,h(t),(f|_{J_{b_1}})^{-1}\circ
h(t),(f|_{J_{b_2}})^{-1}\circ(f|_{J_{b_1}})^{-1}\circ h(t),\ldots)$.
\end{proof} 

\begin{lemma}\label{L1.10}
If $f\in\M$ and $\x\in\invlim f\backslash\{\b\}$
then there is a matchbox $M\subset\invlim f\backslash\{\b\}$ such that
$\x\in\dot{M}$.
\end{lemma}

\begin{proof}
Choose $f\in\M_{n,n}$ to be strictly piecewise monotone and let
$k$ be such that $\pi_k(\x)\ne b$. Let $J$ be a closed arc in $X_n\backslash\{b\}$ 
with $\pi_k(\x)\in\dot{J}$. Then $M=\pi_k^{-1}(J)$ is a matchbox with
$\x\in\dot{M}$.
\end{proof} 

We will say $b$ is {\em distinguished} for $f\in\M$ if $f$ folds at $b$ or
$\#\R(f)\ne 2$. Let $\M^*=\{f\in\M:b\;{\rm is\;distinguished\;for}\; f\}$.

\begin{prop}\label{P1.11}
If $f,g\in\M^*$ and $\psi$ is a homeomorphism from $\invlim f$
to $\invlim g$ then $\psi(\b)=\b$.
\end{prop}

\begin{proof}
This follows immediately from Lemmas \ref{L1.6}-\ref{L1.8} and \ref{L1.10}:
$\b$ is the only point which does not admit a matchbox neighborhood.
\end{proof} 

Let $f\in\M_{n,n}$ be such that $b$ is distinguished for $f$. For each
$i=1,\ldots,n$, let $I_i$ be a closed interval in $\dot{S}_i^1$ and
$h_i:[-1,1]\to I_i$ be a homeomorphism such that the orientation induced on
$I_i$ by the canonical orientation of $[-1,1]$ agrees with the orientation of
${S}_i^1$. Let $\chi_f:\A_n\to\W(\A_n)$ be the substitution determined by
$f$. If $\chi_f(a_i)=c_1\ldots c_s$ then for each $j=0,\ldots,s$, we choose
$x_{i,j}\in[-1,1]$ such that $x_{i,0}=-1$, $x_{i,s}=1$, and
$x_{i,j-1}<x_{i,j}$ for each $j=1,\ldots,s$. Let $\tilde{f}\in\M_{n,n}$ be such 
that $\tilde{f}(X_n\backslash\bigcup_{i=1}^n\dot{I}_i)=\{b\}$ and
$\tilde{f}$ takes $\dot{I}_{i,j}\equiv h_i((x_{i,j-1},x_{i,j}))$
homeomorphically onto $\dot{S}_{\ell}^1$ in an orientation preserving (resp.
reversing) way if $c_j=a_{\ell}$ (resp. $\bar{a}_{\ell}$).
Since $\tilde{f}$ follows the same pattern as $f$ and $b$ is
distinguished for $\tilde{f}$, $\invlim \tilde{f}\backslash\{\b\}\simeq 
\invlim f\backslash\{\b\}$.

\begin{lemma}\label{L1.12}
If $M$ is a matchbox in $\invlim \tilde{f} \backslash\{\b\}$
then there exists a $K$ such that $M\subset\bigcup_{i=1}^n\pi_k^{-1}(I_i)$
for every $k\ge K$.
\end{lemma}

\begin{proof}
For each $\x\in\invlim\tilde{f}\backslash\{\b\}$, so in particular for each
$\x\in M$, there is a $k$ such that $\pi_j(\x)\in\bigcup_{i=1}^n\dot{I}_i$
for all $j\ge k$. Then, by compactness of $M$ there is a $K$ such that 
$\pi_k(\x)\in\bigcup_{i=1}^n\dot{I}_i$ for all $\x\in M$ and all $k\ge K$. 
That is, $M\subset\bigcup_{i=1}^n\pi_k^{-1}(I_i)$ for all $k\ge K$.
\end{proof}

\begin{lemma}\label{L1.13}
If $M$ is an oriented matchbox in $\invlim 
\tilde{f}\backslash\{\b\}$ then there is a positive integer $K$ such that 
for each $k\ge K$ there is a collection of pairwise disjoint matchboxes
$M_1,\ldots,M_{\ell_k}$ whose
union is $M$ such that, for $j=1,\ldots,\ell_k$, the following conditions
hold:
\begin{itemize}
\item[$i)$]there exists an $i\in\{1,\ldots,n\}$ such that 
$M_j\subset\pi_k^{-1}(I_i)$;
\item[$ii)$]$M_j$, with the orientation inherited from $M$, and
$\pi_k^{-1}(I_i)$ are oriented coherently or anticoherently;
\item[$iii)$]there is an $\x_j\in \dot{M}_j$ such that each match of $M_j$ passes
through $\pi_k^{-1}\circ\pi_k(\{\x_j\})$.
\end{itemize}
\end{lemma}

\begin{proof}
Let $h:S\times[-1,1]\to M$ be a parametrization of $M$. By Lemma \ref{L1.12},
there exists an  
$L$ such that $M\subset\bigcup_{i=1}^n\pi_{\ell}^{-1}(I_i)$ for all $\ell\ge L$.
Given $\x\in\dot{M}$ and an open set $U$, where $\x\in U\subset\dot{M}$,
there exists an $\ell\ge L$ such that 
$\pi_{\ell}^{-1}\circ\pi_{\ell}(\{\x\})\subset U$. The matchbox $M_{\x,\ell} 
\equiv h(pr_1\circ h^{-1}\circ\pi_{\ell}^{-1}\circ\pi_{\ell}(\{\x\})
\times[-1,1])\subset M$ inherits an orientation from $M$. For
$\z\in\pi_{\ell}^{-1}\circ\pi_{\ell}(\{\x\})$ and $k\ge\ell$, $M_{\z,k}\subset
M_{\z,\ell}=M_{\x,\ell}$. Also for
$\w,\z\in\pi_{\ell}^{-1}\circ\pi_{\ell}(\{\x\})$ and $k\ge\ell$, $M_{\w,k}$ and
$M_{\z,k}$ are either disjoint or equal. So given $k\ge\ell$, there exist 
$\z_1,\ldots,\z_{r_k}\in\pi_{\ell}^{-1}\circ\pi_{\ell}(\{\x\})$ such that
$M_{\z_1,k},\ldots,M_{\z_{r_k},k}$ is a collection of pairwise disjoint
matchboxes whose union is $M_{\x,\ell}$.
For each $\x\in C\equiv
h(S\times\{0\})$ (cross-section of $M$), choose $\ell(\x)\ge L$ such that
$\pi_{\ell(\x)}^{-1}\circ\pi_{\ell(\x)}(\{\x\})\subset\dot{M}$. So there exist
$\x_1,\ldots,\x_s\in C$ such that
$M_{\x_1,\ell(\x_1)},\ldots,M_{\x_s,\ell(\x_s)}$ is a minimal cover of $M$
(i.e., $M_{\x_i,\ell(\x_i)}\not\subset\bigcup_{j\ne i}M_{\x_j,\ell(\x_j)}$ for each
$i=1,\ldots,s$). Let $M_i$ be the matchbox $M_{\x_i,\ell(\x_i)}\backslash
\bigcup_{j<i}M_{\x_j,\ell(\x_j)}$ and 
$C_i=\{\z\in\pi_{\ell(\x_i)}^{-1}\circ\pi_{\ell(\x_i)}(\{\x_i\}):\z\in M_i\}$.
Choose $k_i\ge\ell(\x_i)$ such that, for all $\z\in C_i$, $M_{\z,k_i}\subset M_i$ 
and $M_{\z,k_i}$ and $\pi_{k_i}^{-1}(I_j)$ are oriented coherently or
anticoherently where $j$ is such that $I_j$ contains $\pi_{k_i}(\z)$. Let
$k\ge K\equiv\max\{k_1,\ldots,k_s\}$. Then there exist
$\z_1,\ldots,\z_{r_k}\in\bigcup_{i=1}^s C_i$ such that 
$M_{\z_1,k},\ldots,M_{\z_{r_k},k}$ is a collection of pairwise disjoint
matchboxes whose union is $M$.
\end{proof} 

\begin{lemma}\label{L1.14}
If $M$is a matchbox in $\invlim \tilde{f}\backslash\{\b\}$ that
intersects each arc component of $\invlim \tilde{f}\backslash\{\b\}$ then
there is a $K$ such that $\pi_k^{-1}(I_i)\cap M\ne\emptyset$ for $k\ge K$ and
$i=1,\ldots,n$.
\end{lemma}

\begin{proof}
Choose $L$ such that $\pi_{\ell}^{-1}\circ\pi_{\ell}(\{\x\})\subset
M\subset\bigcup_{i=1}^n\pi_{\ell}^{-1}(I_i)$ for each $\x$ in a 
cross-section of $M$ and $\ell\ge L$.
For purposes of contradiction, assume there exist $L\le k_1<k_2<\ldots$ and
$I_{i_{k_1}},I_{i_{k_2}},\ldots$ such that $\pi_{k_j}^{-1}(I_{i_{k_j}})\cap
M=\emptyset$. If $\pi_k^{-1}(I_i)\cap M=\emptyset$ for some $i$ and $k\ge L$
then $\pi_{\ell}^{-1}(I_j)\cap M=\emptyset$ for each $I_j$ covered by
$\tilde{f}^{k-\ell}(I_i)$ where $L\le \ell\le k$. So there exist
$i_L,i_{L+1},\ldots$ such that $\pi_j^{-1}(I_{i_j})\cap M=\emptyset$
and $I_{i_j}\subset\tilde{f}(I_{i_{j+1}})$ for $j=L,L+1,\ldots$ . There exists
an   
$\x\in\invlim\tilde{f}\backslash\{\b\}$ such that $\pi_j(\x)\in I_{i_j}$ for
$j=L,L+1,\ldots$ . Choose $\z\in M$ so that $\x$ and $\z$ lie on the same arc
component of $\invlim\tilde{f}\backslash\{\b\}$ and let $A$ be the arc between
$\x$ and $\z$. There is then an $\ell\ge L$ such that $b\not\in\pi_{\ell}(A)$. 
Thus $\pi_{\ell+1}(A)\subset I_{i_{\ell+1}}$. But this contradicts
$\pi_{\ell+1}^{-1}(I_{i_{\ell+1}})\cap M=\emptyset$.
\end{proof} 

Let $g\in\M_{m,m}$ and define the map $\tilde{g}:X_m\to X_m$ as above.

\begin{lemma}\label{L1.15}
If $\invlim \tilde{f}\backslash\{\b\}\simeq\invlim \tilde{g}\backslash\{\b\}$
then $\chi_{\tilde{f}}\sim_w\chi_{\tilde{g}}$.
\end{lemma}
\begin{proof}
Let $\Psi$ be a homeomorphism from $\invlim
\tilde{f}\backslash\{\b\}$ to $\invlim \tilde{g}\backslash\{\b\}$. We inductively
choose integers $0=s_1<t_1<s_2<t_2<\ldots$ in the following way. For a given  
$s_k$, $\pi_{s_k}^{-1}(I_1),\ldots,\pi_{s_k}^{-1}(I_n)$ are pairwise disjoint
matchboxes in $\invlim \tilde{f}\backslash\{\b\}$, with orientation given by
Lemma \ref{L1.9},
and $\bigcup_{\ell=1}^n\pi_{s_k}^{-1}(I_{\ell})$ intersects each arc component of
$\invlim \tilde{f}\backslash\{\b\}$. So $\Psi(\pi_{s_k}^{-1}(I_1)),\dots,
\Psi(\pi_{s_k}^{-1}(I_n))$ are pairwise disjoint oriented matchboxes in $\invlim
\tilde{g}\backslash\{\b\}$ and $\Psi(\bigcup_{\ell=1}^n\pi_{s_k}^{-1}(I_{\ell}))=
\bigcup_{\ell=1}^n\Psi(\pi_{s_k}^{-1}(I_{\ell}))$ intersects each arc component of
$\invlim \tilde{g}\backslash\{\b\}$. By Lemmas \ref{L1.13} and \ref{L1.14}, 
we can choose a $t_k>s_k$ such that, for each $i=1,\ldots,m$,
$\pi_{t_k}^{-1}(I_i^{\prime})\cap\Psi(\bigcup_{\ell=1}^n\pi_{s_k}^{-1}(I_{\ell}))
\ne\emptyset$ and, for each $j=1,\ldots,n$, $\Psi(\pi_{s_k}^{-1}(I_j))$ is
the union of a finite collection of pairwise disjoint matchboxes satisfying 
conditions $i)-iii)$ of Lemma \ref{L1.13}. Similarly, for a given $t_k$,
$\Psi^{-1}(\pi_{t_k}^{-1}(I_1^{\prime})),\ldots,
\Psi^{-1}(\pi_{t_k}^{-1}(I_m^{\prime}))$ are pairwise disjoint oriented matchboxes 
in $\invlim \tilde{f}\backslash\{\b\}$ and
$\Psi^{-1}(\bigcup_{\ell=1}^m\pi_{t_k}^{-1}(I_{\ell}^{\prime}))$ intersects each
arc component of $\invlim \tilde{f}\backslash\{\b\}$. By Lemmas \ref{L1.13}
and \ref{L1.14}, we can choose an $s_{k+1}>t_k$
such that, for $j=1,\ldots,n$, $\pi_{s_{k+1}}^{-1}(I_j)\cap\Psi^{-1}
(\bigcup_{\ell=1}^m\pi_{t_k}^{-1}(I_{\ell}^{\prime}))\ne\emptyset$ and, for
each $i=1,\ldots,m$, $\Psi^{-1}(\pi_{t_k}^{-1}(I_i^{\prime}))$ is the union
of a finite collection of pairwise disjoint matchboxes satisfying conditions
$i)-iii)$ of Lemma \ref{L1.13}. For $k=1,2\ldots,$ we define a substitution
$\sigma_k:\A_m\to\W(\A_n)$ in the following way. Choose a match
$A_i^{\prime}$ of $\pi_{t_k}^{-1}(I_i^{\prime})$. Let
$A_{i,1}^{\prime}<\ldots<A_{i,s}^{\prime}$ be the arc components of
$A_i^{\prime}\cap\Psi(\bigcup_{\ell=1}^n\pi_{s_k}^{-1}(I_{\ell}))\ne\emptyset$,
with the ordering given by the positive orientation of
$\pi_{t_k}^{-1}(I_i^{\prime})$.
Let $\sigma_k(a_i^{\prime})=c_1\ldots c_s$ where $c_j=a_{\ell}$ (resp.
$\bar{a}_{\ell}$) if
$\Psi^{-1}(A_{i,j}^{\prime})\subset\pi_{s_k}^{-1}(I_{\ell})$ and 
$\Psi^{-1}(A_{i,j}^{\prime})$ and $\pi_{s_k}^{-1}(I_{\ell})$ are oriented
coherently (resp. anticoherently). Let $\sigma_k(\bar{a}_i^{\prime})=
\overline{\sigma_k(a_i^{\prime})}$. Conditions $i)-iii)$ of Lemma
\ref{L1.13} ensure that this substitution is independent of the chosen
match. We define
$\tau_k:\A_n\to\W(\A_m)$ similarly. Let $n_k=s_{k+1}-s_k$ and
$m_k=t_{k+1}-t_k$. If $A$ is a match of
$\bigcup_{\ell=1}^n\pi_{s_{k+1}}^{-1}(I_{\ell})$, then the substitution
$\sigma_k\circ\tau_k:\A_n\to\W(\A_n)$ describes the way in which
$\Psi^{-1}\circ\Psi(A)$ passes through
$\bigcup_{\ell=1}^n\pi_{s_k}^{-1}(I_{\ell})$, but
$(\chi_{\tilde{f}})^{n_k}=(\chi_{\tilde{f}})^{s_{k+1}-s_k}:\A_n\to\W(\A_n)$
also describes the way in which $A$ passes through
$\bigcup_{\ell=1}^n\pi_{s_k}^{-1}(I_{\ell})$. So
$\sigma_k\circ\tau_k=(\chi_{\tilde{f}})^{n_k}$. Similarly
$\tau_k\circ\sigma_{k+1}=(\chi_{\tilde{g}})^{m_k}$. Thus 
$\chi_{\tilde{f}}\sim_w\chi_{\tilde{g}}$.
\end{proof} 

\begin{theorem}\label{T1.16}
Given $f,g\in\M^*$, $\invlim f\simeq\invlim g$ iff
$\chi_f\sim_w\chi_g$.
\end{theorem}

\begin{proof}
By Theorem \ref{T1.5}, if $\chi_f\sim_w\chi_g$ then $\invlim
f\simeq\invlim g$. Conversely, if $f,g\in\M^*$ and $\invlim f\simeq\invlim g$,
then $\invlim f\backslash\{\b\}\simeq\invlim g\backslash\{\b\}$ by Proposition
\ref{P1.11}. Let $\tilde{f}$, $\tilde{g}$ be as above.
Then $\invlim\tilde{f}\backslash\{\b\}\simeq\invlim
f\backslash\{\b\}$ and $\invlim g\backslash\{\b\}\simeq
\invlim\tilde{g}\backslash\{\b\}$. So by Lemma \ref{L1.15}, 
$\chi_f=\chi_{\tilde{f}}\sim_w\chi_{\tilde{g}}=\chi_g$.
\end{proof} 
  
Let $\M^I=\M\backslash\M^*$. The $I$ in the notation stands for immersion:
$X_n$ can be given the structure of a smooth branched manifold and for each
$f\in\M^I$ there is a smooth immersion $g\in\M^I$ with $\chi_f=\chi_g$. But
this will not concern us here (see \cite{W2},\cite{J}). We will say that
$f\in\M_n^I\equiv\M_{n,n}\cap\M^I$ is {\em orientation preserving}
provided, with respect to some orientation on the $S_i^1$,
$\chi_f(a_i)\in\W(\{a_1,\ldots,a_n\})$ for each $i\in\{1,\ldots,n\}$. Let
$\M_o^I=\{f\in\M^I:f^2\;{\rm is\;orientation\;preserving}\}$. 
The following theorem is proved in \cite{J}.

\begin{theorem}\label{T1.17}
If $f,g\in\M_o^I$ and $\invlim f\simeq\invlim g$ then $A_f\sim_wA_g$.
\end{theorem}

It is also proved in \cite{J} that for $f\in\M^I\backslash\M_o^I$, 
say $f:X_n\to X_n$,
there is an associated $f^{\prime}\in\M_o^I$, $f^{\prime}:X_{m}\to
X_{m}$. Moreover, if $\invlim f\simeq\invlim g$ for 
$f,g\in\M^I\backslash\M_o^I$, 
then $\invlim f^{\prime}\simeq\invlim g^{\prime}$, which with Theorem 
\ref{T1.17}, yields the following.

\begin{theorem}\label{T1.18}
If $f,g\in\M^I\backslash\M_o^I$ and $\invlim f\simeq\invlim g$ then
$A_{f^{\prime}}\sim_wA_{g^{\prime}}$.
\end{theorem}

A square matrix $A$ with nonnegative entries is {\em aperiodic} provided
$A^k>0$ for some $k$. Such a matrix has a simple real eigenvalue $\lambda_A$,
greater in modulus than all other eigenvalues, called the {\em Perron eigenvalue}
of $A$. Let $\mathbb{Q}(\lambda_A)$ denote the extension field of the rationals
$\mathbb{Q}$ by $\lambda_A$.

\begin{lemma}\label{L1.19}
If $A,B$ are aperiodic square nonnegative integral matrices and
$A\sim_wB$ then $\mathbb{Q}(\lambda_A)=\mathbb{Q}(\lambda_B)$.
\end{lemma}

\begin{proof}See \cite{BD}.\end{proof}

\begin{corollary}\label{C1.20}
Let $f,g\in\M$ and suppose that $A_f$ and $A_g$ are aperiodic. If $\invlim
f\simeq\invlim g$ then $\mathbb{Q}(\lambda_A)=\mathbb{Q}(\lambda_B)$.
\end{corollary}

\begin{proof}
If $f,g\in\M^*$, this follows from Lemmas \ref{L1.1}, \ref{L1.19} and 
Theorem \ref{T1.16}. If $f,g\in\M_o^I$, this follows from Theorem
\ref{T1.17} and Lemma \ref{L1.19}. If $f,g\in\M^I\backslash\M_o^I$, 
Jacklitch \cite{J}
proves that $\lambda_{f^{\prime}}=\lambda_f$, $\lambda_{g^{\prime}}=\lambda_g$,
so the corollary follows from Theorem \ref{T1.18} and Lemma \ref{L1.19}. There
are no other possibilities for $\invlim f\simeq\invlim g$.
\end{proof} 

The proofs of Theorems \ref{T1.17} and \ref{T1.18} are quite different.
Nonetheless, we make the following conjecture.

\begin{conj}\label{C1.21}
Given $f,g\in\M$, $\invlim f\simeq\invlim g$ iff $\chi_f\sim_w\chi_g$.
\end{conj}

\section{Applications}
In this section we present several applications of the preceding theorems.
The first of these will be to one-dimensional substitution tilings and for
this we borrow heavily from \cite{AP}.

Consider a substitution $\chi:\A=\{a_1,\ldots,a_n\}\to\W(\A)$, $n\ge 2$, with
associated matrix $A_{\chi}$. To avoid degeneracy, we assume that for no word
$w$ and positive integer $k$ is $\chi^k(w)=ww\ldots w$.
We also assume that $A_{\chi}$ is aperiodic
with Perron eigenvalue $\lambda_{\chi}$ and associated positive (left)
eigenvector $\vec{v}_{\chi}$ with entries $\lambda_1,\ldots,\lambda_n$. The
intervals $P_i=[0,\lambda_i]$, $i=1,\ldots,n$, are called prototiles (consider
$P_i$ to be distinct from $P_j$ for $i\ne j$ even if $\lambda_i=\lambda_j$).
A {\em tiling} $T$ of $\mathbb{R}$ by the prototiles is a collection $T=
\{T_i\}_{i=-\infty}^{\infty}$ of tiles $T_i$,
$\bigcup_{i=-\infty}^{\infty}T_i=\mathbb{R}$, with each $T_i$ a translate of
some $P_j$ and with $T_i\cap T_{i+1}$ a singleton for each $i$.

Now $\lambda_{\chi}\lambda_i=\sum_{j=1}^r\lambda_{i(j)}$ where
$i(j)=k$ provided $\chi(a_i)=c_1\ldots c_r$ with $c_j=a_k$. Thus
$\lambda_{\chi}P_i=\bigcup_{s=0}^{r-1}(\sum_{j=1}^s\lambda_{i(j)})+P_{i(s+1)}$ 
that is,
$\lambda_{\chi}P_i$ is tiled by $\{T_s\}_{s=0}^{r-1}$,
$T_s=(\sum_{j=1}^s\lambda_{i(j)})+P_{i(s+1)}$. This process is called inflation and
substitution and extends to a map $\sigma_{\chi}$ taking a tiling
$T=\{T_i\}_{i=-\infty}^{\infty}$ of $\mathbb{R}$ by prototiles to a new tiling,
$\sigma_{\chi}(T)$ of $\mathbb{R}$ by prototiles defined by inflating, substituting, 
and suitably translating each $T_i$.

There is a natural topology on the collection $\Sigma_{\chi}$ of all tilings
of $\mathbb{R}$ by prototiles ($\{T_i\}_{i=-\infty}^{\infty}$ and
$\{T_i^{\prime}\}_{i=-\infty}^{\infty}$ are ``close" if there is an
$\epsilon$ near zero and a large $N$ so that if $0\in T_k\cap
T_{\ell}^{\prime}$ then $T_{k+i}=\epsilon+ T_{\ell+i}^{\prime}$ for
$i=-N,\ldots,N$, see \cite{AP} for details). 
$\Sigma_{\chi}$ is a compact metrizable space with this
topology and $\sigma_{\chi}:\Sigma_{\chi}\to\Sigma_{\chi}$ is continuous. Let
$T\in\Sigma_{\chi}$ be any tiling that is periodic under $\sigma_{\chi}$
(there are infinitely many such) and consider the orbit, $\O(T)$, of $T$
under the flow $\phi_t:\Sigma_{\chi}\to\Sigma_{\chi}$ by
$\phi_t(\{T_i\}_{i=-\infty}^{\infty})=\{t+T_i\}_{i=-\infty}^{\infty}$. Let
$\T_{\chi}=cl(\O(T))$; $\T_{\chi}$ is the {\em tiling space associated with}
$\chi$. Then $\T_{\chi}$ does not depend on the $T$ chosen above, $\phi_t$ is
minimal on $\T_{\chi}$ (so in particular, $\T_{\chi}$ is a continuum) and
$\sigma_{\chi}:\T_{\chi}\to\T_{\chi}$ is a homeomorphism.

Given two such substitutions $\chi$, $\psi$, we wish to apply theorems from
the previous section in an effort to distinguish $\T_{\chi}$ from $\T_{\psi}$ if
they are nonhomeomorphic. In order to (partially) accomplish this, we need to
consider {\em collared} matrices $A_{\tilde{\chi}}$ and $A_{\tilde{\psi}}$
associated with $\chi$ and $\psi$. Consider the collection $\tilde{\A}=
\tilde{\A}_{\chi}=\{(a_ia_ja_k): a_ia_ja_k \;{\rm occurs\;as\;a\;subword\;of}
\;\chi^{\ell}(a_s),\;{\rm some}\;\ell,s\}$.
Define $\tilde{\chi}:\tilde{\A}\to\W(\tilde{\A})$ by
$\tilde{\chi}((a_ia_ja_k))=(b_0b_1b_2)(b_1b_2b_3)\ldots
(b_{s-2}b_{s-1}b_s)(b_{s-1}b_sb_{s+1})$
where $\chi(a_i)=\ldots b_0$, $\chi(a_j)=b_1\ldots b_s$, and
$\chi(a_k)=b_{s+1}\ldots$ . Let $A_{\tilde{\chi}}$ be the matrix associated
with $\tilde{\chi}$ (choose any numbering of the elements of $\tilde{\A}$).
Then $A_{\tilde{\chi}}$ is also aperiodic. Let $\lambda_{\tilde{\chi}}$ be the
Perron eigenvalue of $A_{\tilde{\chi}}$.

\begin{lemma}\label{L2.1}
$\lambda_{\tilde{\chi}}=\lambda_{\chi}$.
\end{lemma}

\begin{proof}
By construction of $\tilde{\chi}$, there is a nonnegative integer matrix $P$
of rank $n$ such that $PA_{\tilde{\chi}}=A_{\chi}P$. Let 
$\tilde{v}$ be a positive
eigenvector of $A_{\tilde{\chi}}$ corresponding to the eigenvalue
$\lambda_{\tilde{\chi}}$. Then
$A_{\chi}P\tilde{v}=PA_{\tilde{\chi}}\tilde{v}
=P\lambda_{\tilde{\chi}}\tilde{v}
=\lambda_{\tilde{\chi}}P\tilde{v}$.
Thus $v=P\tilde{v}$ is a positive eigenvector for
$A_{\chi}$ with eigenvalue $\lambda_{\tilde{\chi}}$.
Thus $\lambda_{\tilde{\chi}}=\lambda_{\chi}$
by the Perron-Frobenius Theorem.
\end{proof}

\begin{theorem}\label{T2.2}
If $\T_{\chi}\simeq\T_{\psi}$ then $A_{\tilde{\chi}}\sim_w A_{\tilde{\psi}}$.
\end{theorem}
\begin{proof}
Anderson and Putnam (\cite{AP}) prove that there are branched
one-manifolds $G_{\chi}$ and $G_{\psi}$ and expanding, orientation preserving
immersions $f_{\chi}:G_{\chi}\to G_{\chi}$ and $f_{\psi}:G_{\psi}\to
G_{\psi}$ such that $\T_{\chi}\simeq\invlim f_{\chi}$ and
$\T_{\psi}\simeq\invlim f_{\psi}$. Furthermore, the Markov matrices
associated with $f_{\chi}$ and $f_{\psi}$ are $A_{\tilde{\chi}}$ and
$A_{\tilde{\psi}}$, respectively. Using the Williams moves to split the branched
manifolds into wedges of circles, there are $g_{\chi}:X_n\to X_n$ and
$g_{\psi}:X_m\to X_m$ with $\invlim g_{\chi}\simeq\invlim
f_{\chi}$ and $\invlim g_{\psi}\simeq\invlim f_{\psi}$ (see \cite{W2} and 
\cite{J}). Moreover, $g_{\chi},g_{\psi}\in\M_{o.p.}^I$ and the associated 
matrices $A_{g_{\chi}}$ and $A_{g_{\psi}}$ are shift equivalent to some 
powers of $A_{\tilde{\chi}}$ and $A_{\tilde{\psi}}$, respectively. In particular 
$A_{g_{\chi}}\sim_w A_{\tilde{\chi}}$ and $A_{g_{\psi}}\sim_w
A_{\tilde{\psi}}$. Assuming $\T_{\chi}\simeq\T_{\psi}$, it follows from
Theorem \ref{T1.17} that $A_{\tilde{\chi}}\sim_w A_{\tilde{\psi}}$.
\end{proof}

\begin{corollary}\label{C2.3} 
If $\T_{\chi}\simeq\T_{\psi}$ then $\mathbb{Q}(\lambda_{\chi})=
\mathbb{Q}(\lambda_{\psi})$.
\end{corollary}

\begin{proof}
This follows from Lemmas \ref{L1.19}, \ref{L2.1} and Theorem \ref{T2.2}.
\end{proof} 

Restricting Conjecture \ref{C1.21} to this setting we have

\begin{conj}\label{C2.4} 
$\T_{\chi}\simeq\T_{\psi}$ iff $\tilde{\chi}\sim_w\tilde{\psi}$.
\end{conj}

The next application we consider is to one-dimensional unstable manifolds of
hyperbolic basic sets. We will assume that combinatorial descriptions of the
unstable manifolds are given in terms of collections of Markov rectangles
satisfying certain felicitous properties. This situation arises, for example,
in \cite{BH} and \cite{BL}. Perhaps the properties R1)-R4) below can
always be arranged for the right choice of rectangles in the presence of a
local product structure (\cite{HK}).

Suppose that $F:M\to M$ is a diffeomorphism of the $k$-dimensional manifold $M$
with zero-dimensional hyperbolic invariant set $\Lambda_F$ having one-dimensional
unstable manifold $W^u(\Lambda_F)$. In addition, we suppose that there
are pairwise disjoint {\em rectangles} $R_1,\ldots,R_n\subset M$ with
$\Lambda_F\subset\bigcup_{i=1}^n\dot{R}_i$ having the following properties.
\begin{itemize}
\item[R1)]Each $R_i$ is parametrized by a homeomorphism
$\phi_i:[-1,1]\times[-1,1]^{k-1}\to R_i$.
\item[R2)]$\Lambda_F=\bigcap_{\ell=-\infty}^{\infty}F^{\ell}(\bigcup_{i=1}^nR_i)$
and $R_i\cap\Lambda_F\ne\emptyset$ for each $i$.
\item[R3)]If $z=\phi_i(t,x)\in R_i$ and
$F(z)=\phi_j(t^{\prime},x^{\prime})\in R_j$, then
$F(\phi_i(\{t\}\times[-1,1]^{k-1}))\subset\phi_j(\{t\}^{\prime}\times[-1,1]^{k-1})$.
\item[R4)]For $z\in\Lambda_F$, let $W^u(z,R_i)$ be the component of
$W^u(z)\cap R_i$ containing $z$. Then $W^u(z,R_i)$ meets each
$\phi_i(\{t\}\times[-1,1]^{k-1})$ exactly once and if $F(z)\in R_j$ then
$F(W^u(z,R_i))\supset W^u(F(z),R_j)$.
\item[R5)]If
$z\in(\bigcup_{i=1}^nR_i)\cap(F^{-1}(M\backslash(\bigcup_{i=1}^nR_i))$
then $F^{\ell}(z)\in M\backslash(\bigcup_{i=1}^nR_i)$ for all $\ell\ge 1$.
\end{itemize}
Given $(M,F,\Lambda_F,\{R_i\}_{i=1}^n)$ satisfying R1)-R5), let $P_i:R_i\to
L_i=[-1,1]\times\{i\}$ be defined by $P_i(\phi_i(t,x))=(t,i)$. For each 
$i\in\{1,\ldots,n\}$,
let $J_1^i<J_2^i<\ldots<J_{k_i}^i$ be the components of $P_i(R_i\cap 
F^{-1}(\bigcup_{j=1}^nR_j))$ (ordered according to the usual order on
$[-1,1]$). Let $f:\bigcup_{i=1}^n(\bigcup_{j=1}^{k_i}J_j^i)\to
\bigcup_{i=1}^nL_i$ be defined by $f((t,i))=P_j(F(z))$ where $P_i(z)=(t,i)$ and 
$F(z)\in R_j$. Then $f$ is
a continuous surjection. Moreover, for each $j\in\{1,\ldots,k_i\}$ and
$i\in\{1,\ldots,n\}$ there is an $L_{\ell}$ such that $f|_{J_j^i}:J_j^i\to
L_{\ell}$ is a continuous monotone surjection. We define the substitution
$\chi_F:\A_n\to\W(\A_n)$ by $\chi_F(a_i)=b_1\ldots b_{k_i}$ with
$b_j=a_{\ell}$ if $f(J_j^i)=L_{\ell}$ and $f$ is orientation preserving 
on $J_j^i$ and $b_j=\bar{a}_{\ell}$ if $f(J_j^i)=L_{\ell}$
and $f$ is orientation reversing on $J_j^i$.

\begin{prop}\label{P2.5}
If $(M,F,\Lambda_F,\{R_i\}_{i=1}^n)$ satisfies R1)-R5) and 
$f^{\prime}\in\M_{n,n}$ follows the pattern of $\chi_F$ then $\invlim
f^{\prime}\backslash\{\b\}\simeq W^u(\Lambda_F)$.
\end{prop}

\begin{proof}
Let $h_i:L_i\to I_i\subset\dot{S}_i^1\subset X_n$ be 
orientation preserving homeomorphisms for $i=1,\ldots,n$. Let
$f^{\prime}:\bigcup_{i=1}^nh_i(\bigcup_{j=1}^{k_i}J_j^i)\to\bigcup_{i=1}^nI_i$
by $f^{\prime}(t)=h_j\circ f\circ h_i^{-1}(t)$, $i$ and $j$ as appropriate.
Extend $f^{\prime}$ to a map in $\M_{n,n}$ that follows the pattern of $\chi_F$
and so that for $x\not\in\I\equiv\bigcup_{i=1}^nI_i$ there is an $\ell$ such
that $(f^{\prime})^{\ell}(x)=b$.
Let $\R=\bigcup_{i=1}^nR_i$ and $P:\R\to\I$ be defined by $P(z)=h_i\circ P_i(z)$.
Then the diagram\\
\begin{center}
\begin{picture}(2.9,2.2)
\put(.4,.4){\makebox(0,0){$(f^{\prime})^{-1}(\I)\cap\I$}}
\put(.4,1.8){\makebox(0,0){$F^{-1}(\R)\cap\R$}}
\put(2.5,.4){\makebox(0,0){$\I$}}
\put(2.5,1.8){\makebox(0,0){$\R$}}
\put(.2,1.1){\makebox(0,0){$P$}}
\put(2.7,1.1){\makebox(0,0){$P$}}
\put(1.8,.2){\makebox(0,0){$f^{\prime}$}}
\put(1.8,2){\makebox(0,0){$F$}}
\put(.4,1.6){\vector(0,-1){1}}
\put(2.5,1.6){\vector(0,-1){1}}
\put(1.3,.4){\vector(1,0){1}}
\put(1.3,1.8){\vector(1,0){1}}
\end{picture}
\end{center}
commutes. Given $z\in W^u(\Lambda_F,\R)\equiv
\bigcup_{i=1}^n(\bigcup_{z\in\Lambda_F\cap R_i}(W^u(z,R_i)))$, let
$H(z)=(P(z),P(F^{-1}(z)),P(F^{-2}(z)),\ldots)$. Then $H(z)\in\invlim
f^{\prime}$ and $\pi_{\ell}(H(z))\in\I$ for $\ell=0,1,\ldots$ .

\noindent{\em Claim:} $H:W^u(\Lambda_F,\R)\to\pi_0^{-1}(\I)$ is a
homeomorphism.

\noindent{\em Proof:} Suppose that $H(z)=H(w)$ for some
$z,w\in W^u(\Lambda_F,\R)$. Then $P(F^{-\ell}(z))=P(F^{-\ell}(w))$ for
all $\ell$. Then $F^{-t}(z),F^{-t}(w)\in\bigcap_{\ell\ge
0}F^{\ell}(P^{-1}(P(F^{-t-\ell}(\{z\})))$ for each $t=0,1,\ldots$ . Given
$\epsilon>0$ there is an $N$ such that diam$(F^N(W^s(v,R_i)))<\epsilon$
for all $v\in\Lambda_F\cap R_i$ and $i=1,\ldots,n$. So there is a $\delta>0$
such that if the Hausdorff distance between $P^{-1}(P(F^{-N}(\{z\})))$ and
$W^s(v,R_i)$ is less than $\delta$ then 
diam$(F^N(P^{-1}(P(F^{-N}(\{z\}))))<\epsilon$. Since
there are $z^{\prime},w^{\prime}\in\Lambda_F$ with $z\in W^u(z^{\prime})$
and $w\in W^u(w^{\prime})$, for sufficiently large $\ell$, the Hausdorff
distance between $P^{-1}(P(F^{-\ell}(\{z\})))$ and
$W^s(F^{-\ell}(z^{\prime}),R_{i_{\ell}})$ and between
$P^{-1}(P(F^{-\ell}(\{w\})))$ and $W^s(F^{-\ell}(w^{\prime}),R_{i_{\ell}})$
is less than $\delta$ ($i_{\ell}$ such that $F^{-\ell}(z^{\prime}),
F^{-\ell}(w^{\prime})\in R_{i_{\ell}}$). It follows that the distance between
$F^{-t}(z)$ and $F^{-t}(w)$ goes to zero as $t\to\infty$. Thus
$z^{\prime}=w^{\prime}$. Since $P(z)=P(w)$ and $W^u(z^{\prime},R_i)$
hits each $\{t\}\times[-1,1]^{k-1}=P_i^{-1}(\{t\})$ in a unique point, $z=w$.
$H$ is surjective because $f^{\prime}$ follows the pattern
$\chi_F:\pi_0^{-1}(\I)=\{\x\in\invlim f^{\prime}:\pi_k(\x)\in\I\;\forall k\}$
so given $\x\in\pi_0^{-1}(\I)$ pick $z\in\bigcap_{\ell\ge
0}F^{-\ell}(P^{-1}(\{x_{\ell}\}))$; then $H(z)=\x$. Note that for $z\in
H^{-1}(W^u(\Lambda_F,\R))$, $H\circ F(z)=\hat{f}^{\prime}(H(z))$. Extend
$H$ to all of $W^u(\Lambda_F)$ by conjugation; that is, for
$z\in W^u(\Lambda_F)$ let $\ell$ be such that
$F^{-\ell}(z)\in W^u(\Lambda_F,\R)$. Put
$H(z)=(\hat{f}^{\prime})^{\ell}(H(F^{-\ell}(z))$. Then (by the above note)
$H$ is well defined on $W^u(\Lambda_F)$ and, by assumption R5), $H$ is
continuous. Clearly $H(W^u(\Lambda_F))=\bigcup_{\ell\ge
0}(\hat{f}^{\prime})^{\ell}(\pi_0^{-1}(\I))=\invlim
f^{\prime}\backslash\{\b\}$. If $\{\z^{\ell}\}_{\ell=1}^{\infty}$ is a
sequence in $W^u(\Lambda_F)$ without a limit point in
$W^u(\Lambda_F)$ then the (minimal) $j_{\ell}$ for which $\z^{\ell}\in
F^{j_{\ell}}(W^u(\Lambda_F,\R))$ diverge to $\infty$. Thus the sequence
$H(\z^{\ell})\to\b$ as $\ell\to\infty$ and it follows $H^{-1}$ is continuous.
\end{proof}

With the assumptions of the proposition, we see that the one-point
compactification $(W^u(\Lambda_F))^*$ is homeomorphic with $\invlim
f^{\prime}$. We will say that $W^u(\Lambda_F)$ is {\em distinguished at}
$\infty$ if $f^{\prime}\in\M^*$.

For $z\in\Lambda_F$, an arc component $A$ of $W^u(z)\backslash\{z\}$ is
a {\em free separatrix} of $W^u(\Lambda_F)$ provided
$A\cap\Lambda_F=\emptyset$. The reader may check that the number of free
separatrices of $W^u(\Lambda_F)$ is precisely $\#\R(f^{\prime})$,
the cardinality of the eventual range of $f_*^{\prime}$ acting on 
the edge germs of $X_n$. Thus we see
that, for $(M,F,\Lambda_F,\{R_i\}_{i=1}^n)$ satisfying R1)-R5),
$W^u(\Lambda_F)$ is distinguished at $\infty$ if and only if
$W^u(\Lambda_F)$ has other than exactly two free separatrices or, for
some $\ell,i,j$, $\chi_F^{\ell}(a_i)$ contains consecutive letters $a_j$ and
$\bar{a}_j$ (that is, $F^{\ell}(R_i)$ passes through $R_j$ then folds back
immediately through $R_j$).

\begin{theorem}\label{T2.6}
Suppose that $(M,F,\Lambda_F,\{R_i\}_{i=1}^n)$ and
$(N,G,\Lambda_G,\{R_i^{\prime}\}_{i=1}^m)$ satisfy R1)-R5) and in addition
suppose that $W^u(\Lambda_F)$ and $W^u(\Lambda_G)$ are
distinguished at $\infty$ then $W^u(\Lambda_F)\simeq W^u(\Lambda_G)$
if and only if $\chi_F\sim_w\chi_G$.
\end{theorem}

\begin{proof}
If $W^u(\Lambda_F)\simeq W^u(\Lambda_G)$ then
$(W^u(\Lambda_F))^*\simeq(W^u(\Lambda_G))^*$ so that, with
$f^{\prime}$, $g^{\prime}$ as in Proposition \ref{P2.5}, $\invlim
f^{\prime}\simeq\invlim g^{\prime}$. From Theorem \ref{T1.16},
$\chi_{f^{\prime}}\sim_w\chi_{g^{\prime}}$ so that $\chi_F\sim_w\chi_G$.
Conversely, if $\chi_F\sim_w\chi_G$ then
$\chi_{f^{\prime}}\sim_w\chi_{g^{\prime}}$ so that $\invlim
f^{\prime}\simeq\invlim g^{\prime}$, by Theorem \ref{T1.16}. Since
$f^{\prime},g^{\prime}\in\M^*$, a homeomorphism of $\invlim f^{\prime}$ to
$\invlim g^{\prime}$ must take $\b$ to $\b$ (Proposition \ref{P1.11}). Thus
$W^u(\Lambda_F)\simeq\invlim f^{\prime}\backslash\{\b\}\simeq\invlim
g^{\prime}\backslash\{\b\}\simeq W^u(\Lambda_G)$.
\end{proof}

The idea of studying one-dimensional unstable manifolds by means of their
underlying substitutions comes from \cite{V}, in which Vago gives a complete
and easy-to-compute classification of the unstable manifolds of one-rectangle
systems on surfaces, up to homeomorphism and up to conjugacy of the
underlying dynamics.

As a final application, consider the {\em tent maps} $T_s:[0,1]\to[0,1]$ 
defined by
$$T_s(x)=\left\{\begin{array}{ccc}
sx+2-s&,&0\le x\le\frac{s-1}{s}\\
-sx+s&,&\frac{s-1}{s}\le x\le 1
\end{array}\right.,\;{\rm for}\;1<s\le 2.$$
Ingram has conjectured that $\invlim
T_{s_1}\not\simeq\invlim T_{s_2}$ for $1<s_1<s_2\le 2$.
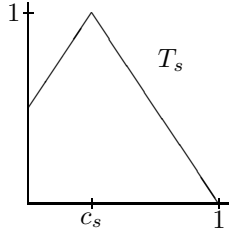
\begin{center}
\begin{picture}(2.4,2.6)
\put(.966,.25){\line(0,1){.1}}
\put(2.298,.25){\line(0,1){.1}}
\put(.25,2.299){\line(1,0){.1}}
\put(.966,.13){\makebox(0,0){$c_s$}}
\put(1.8,1.8){\makebox(0,0){$T_s$}}
\put(2.298,.13){\makebox(0,0){1}}
\put(.15,2.299){\makebox(0,0){1}}
\put(.3,.3){\line(0,1){2.1}}
\put(.3,.3){\line(1,0){2.1}}
\put(.3,1.3){\line(2,3){.666}}
\put(.966,2.299){\line(2,-3){1.332}}
\end{picture}
\end{center}

Consider the set $\P\subset(1,2]$ of parameters for which the critical point
$c_s=\frac{s-1}{s}$ is periodic under $T_s$. Say $s\in\P$ and $c_s$ has
period $n+1$ under $T_s$. Then $\O(c_s)\equiv\{0=x_0<\ldots<x_n=1\}$ breaks
$[0,1]$ into $n$ subintervals $J_i=[x_{i-1},x_i]$, $i=1,\ldots,n$. For such
an $s$, let $\chi_s:\A_n\to\W(\A_n)$ be the substitution defined by:
$\chi_s(a_i)=a_{\ell}a_{\ell+1}\ldots a_{\ell+k}$ if
$T_s(J_i)=J_{\ell}\cup\ldots\cup J_{\ell+k}$ and $T_s$ is increasing on
$J_i$; $\chi_s(a_i)=\bar{a}_{\ell+k}\ldots\bar{a}_{\ell}$ if
$T_s(J_i)=J_{\ell}\cup\ldots\cup J_{\ell+k}$ and $T_s$ is decreasing on
$J_i$.

\begin{theorem}\label{T2.7}
Given $s_1,s_2\in\P$, $\invlim T_{s_1}\simeq\invlim
T_{s_2}$ iff $\chi_{s_1}\sim_w\chi_{s_2}$.
\end{theorem}

\begin{proof}
Given $s\in\P$ with $c_s$ of period $n+1$ under $T_s$, let
$P_s:[0,1]\to X_n=\bigvee_{i=1}^nS_i^1$ take $\dot{J}_i$ (see notation
introduced before this theorem) homeomorphically, and in an orientation 
preserving manner,
onto $\dot{S}_i^1$ for each $i$, and let $P_s(x_i)=b$ for $i=0,\ldots,n$. Then
there is an $f_s\in\M_{n,n}$ such that $P_s\circ T_s=f_s\circ P_s$. Moreover,
$\chi_s=\chi_{f_s}$. Let $\c_s$ denote the point in $\invlim T_s$ with
$\pi_{k(n+1)}(\c_s)=c_s$ for $k=0,1,\ldots$ and let
$\O(\c_s)\equiv\{\hat{T}^k(\c_s):k=0,\ldots,n\}$ be its orbit. Then the map
$\hat{P}_s:\invlim T_s\to\invlim f_s$ defined by
$\hat{P}_s((x_0,x_1,\ldots))=(P_s(x_0),P_s(x_1),\ldots)$ is a surjection and
is one-to-one except on $\O(\c_s)$ where $\hat{P}(\O(\c_s))=\{\b\}$. For
$s\in\P$, the endpoints of the continuum $\invlim T_s$ are exactly the points
of $\O(\c_s)$ (see \cite{BM}). For $s_1,s_2\in\P$, any homeomorphism from $\invlim
T_{s_1}$ to $\invlim T_{s_2}$ must take endpoints to endpoints so such a
homeomorphism would descend to a homeomorphism of $\invlim f_{s_1}$ with
$\invlim f_{s_2}$ by way of $\hat{P}_{s_1}$, $\hat{P}_{s_2}$. Moreover, the
arc components of the endpoints of $\invlim T_s$ are in one-to-one
correspondence with elements of $\R(f_s)$, the eventual range of the action
of $f_s$ on edge germs. Since $c_s$ is not periodic
of period 2 for any $s\in\P$, $\#\R(f_s)\ne 2$ for any $s\in\P$. That
is $f_s\in\M^*$ for $s\in\P$. Thus a homeomorphism of $\invlim f_{s_1}$ with
$\invlim f_{s_2}$ would necessarily take $\b$ to $\b$ and would then lift to
a homeomorphism of $\invlim T_{s_1}$ with $\invlim T_{s_2}$. Thus for
$s_1,s_2\in\P$, $\invlim T_{s_1}\simeq\invlim T_{s_2} \Leftrightarrow \invlim
f_{s_1}\simeq\invlim f_{s_2} \Leftrightarrow \chi_{f_{s_1}}\sim_w\chi_{f_{s_2}}
\Leftrightarrow \chi_{s_1}\sim_w\chi_{s_2}$.
\end{proof} 

The above sort of reasoning can be extended to prove that for $s_1$, $s_2$ in
the {\em post-critically-finite} set ($\O(c_s)$ finite, $1<s\le 2$), $\invlim
T_{s_1}\simeq\invlim T_{s_2}$ if and only if $\chi_{s_1}\sim_w\chi_{s_2}$.
Theorem 2.7 generalizes the result of \cite{BD}: if $s_1,s_2\in\P$ and $\invlim
T_{s_1}\simeq\invlim T_{s_2}$ then $A_{\chi_{s_1}}\sim_w A_{\chi_{s_2}}$.
Bruin \cite{Bru} has recently shown that, for $s_1$ and $s_2$ in the
post-critically finite set, if $\invlim T_{s_1}\simeq\invlim T_{s_2}$ then
$\log s_1/\log s_2 \in\mathbb{Q}$.

\end{document}